\documentclass{amsart}
\usepackage{amssymb}
\usepackage{verbatim}
\usepackage{mathrsfs}
\usepackage{dsfont}
\usepackage{braket}

\newtheorem{theorem}{Theorem}[section]
\newtheorem{lemma}[theorem]{Lemma}

\newtheorem{proposition}[theorem]{Proposition}

\theoremstyle{definition}
\newtheorem{definition}[theorem]{Definition}
\newtheorem{example}[theorem]{Example}
\newtheorem{observation}[theorem]{Observation}

\theoremstyle{remark}
\newtheorem{remark}[theorem]{Remark}

\newcommand{\qA}{{\mathfrak{A}}}
\newcommand{\cB}{{\mathcal B}}

\newcommand{\cH}{{\mathcal H}}
\newcommand{\cL}{{\mathcal L}}
\newcommand{\cK}{{\mathcal K}}

\newcommand{\cS}{{\mathcal S}}

\newcommand{\Rn}{{\rm I\!R}} 
\newcommand{\Nn}{{\rm I\!N}} 
\newcommand{\Cn}{{\setbox0=\hbox{
$\displaystyle\rm C$}\hbox{\hbox
to0pt{\kern0.6\wd0\vrule height0.9\ht0\hss}\box0}}} 

\numberwithin{equation}{section}

\newcommand{\Tr}{\mathrm{Tr}}

\newcommand{\jed}{{\mathbb{I}}}
\begin{document}

\title{On the structure of the set of positive maps.}

\author{W. A. Majewski}

\address{Unit for BMI, North-West-University, Potchefstroom, South Africa.}
\email{fizwam@ug.edu.pl; and fizwam@gmail.com}

\date{\today}




\begin{abstract}
The full description of the set of positive maps $T: \qA \to \cB(\cH)$ ($\qA$ a $C^*$-algebra) is  given.  The approach is based on the simple prescription for selecting various types of positive maps. This prescription stems from the Grothendieck  theory of projective tensor products  complemented by  the theory of tensor connes.
 In particular,
the origin of non-decomposable maps is clarified.
\end{abstract}

\maketitle
\vspace*{1.5cm}
\noindent

\section{Introduction}
Despite the facts that positive maps are an essential ingredient in a description of quantum systems and that they play an  important role in mathematics \cite{paulsen}, \cite{St}, a characterization of the structure of the set of all positive maps has been a long standing challenge in mathematical physics. The key reason behind that is the complexity of this structure - the structure of positive maps is drastically nontrivial even for the finite dimensional case. To illustrate this point it is enough to note that even
the convex structure of the positive maps, $\Phi : \cB(\cH) \to \cB(\cH)$, is highly complicated even in very low dimensions of the Hilbert space $\cH$. 

In the sixties it was shown \cite{St1} (see also St{\o}rmer's book \cite{St} and references given there) that every positive map for the 2D case, i.e. for $\dim \cH = 2$, is decomposable. The first example of non-decomposable map was given by Choi \cite{Choi}, see also \cite{Choi1} and \cite{Wor}, for 3D case, i.e. for $\dim \cH = 3$. Since then, other examples of non-decomposable maps were constructed. In particular, by results of Woronowicz \cite{Wor} and St{\o}rmer \cite{St1}, if $\dim \cH_1 \cdot \dim\cH_2 \leq 6$, all positive maps $T:\cB(\cH_1) \to \cB(\cH_2)$ are decomposable but \textit{this is not true in higher dimensions.} On the other hand, the emergence of non-decomposable maps may be considered as a huge obstacle in getting a canonical form for a positive map.

The present work being a continuation of our previous papers \cite{WAM}, \cite{WAM2}, \cite{LMM}, \cite{WAM3}, \cite{MMO},  \cite{WAM4} and \cite{WAM_o},  provides an analysis  as well as a description of the structure of the set of positive maps.  These results stem from the observation that the linear tensor product structure is not compatible with either topology (there are many cross-norms) or an order (there are many tensor cones).

The deep Grothendieck's result gives the relation between linear maps on Banach spaces $\cL(X,Y)$ and linear continuous functionals on the projective tensor product $(X \otimes_{\pi} Y_*).$ That is a landmark, a crucial ingredient of our approach. An analysis of various orders in $(X \otimes_{\pi} Y_*)$ yields essential information about the structure of positive maps. In particular,
as a byproduct,  we will get an explanation
of the origin of non-decomposable maps. 

This paper is organized as follows: first we give necessary preliminaries in Section 2. Next, in Section 3, a description of various types of positive maps will be given. As a result, the structure of the set of all positive maps will be characterized.  Conclusions and final remarks are given in Section 4.

\section{Definitions and notations}
For any $C^*$-algebra $\mathfrak{A}$, we denote the set of all self-adjoint (positive) elements of $\mathfrak{A}$ by $\mathfrak{A}_h$ ($\mathfrak{A}^+$). 
If $\mathfrak{A}$ is a unital $C^*$-algebra then a state on $\mathfrak{A}$ is a linear functional $\phi: \mathfrak{A} \to \Cn$ such that $\phi(a) \geq 0$ for every $a \geq 0$ ($ a \in \mathfrak{A}^+$) and $\phi(\jed) = 1$, where $\jed$ is the unit of $\mathfrak{A}$. The set of all states on $\mathfrak{A}$ will be denoted by $\cS_{\mathfrak{A}}$. In particular
$$(\mathfrak{A}_h, \mathfrak{A}^+)$$
is an ordered Banach space.

A linear map $T: \qA_1 \to \qA_2$ between $C^*$-algebras $\qA_1$ and $\qA_2$ is called positive if $T(\qA_1^+) \subseteq \qA_2^+$. For $k\in \Nn$ we consider a map $T_k: M_k(\Cn) \otimes \qA_1 \to M_k(\Cn) \otimes \qA_2$ where $M_k(\Cn)$ denotes the algebra of $k\times k$ matrices with complex entries and $T_k = id_{M_k} \otimes T$. We say that $T$ is k-positive if the map $T_k$ is positive. The map $T$ is said to be completely positive (cp for short) if $T$ is k-positive for every $k \in \Nn$. The set of all linear, bounded (unital) positive maps $T: \mathfrak{A}_1 \to \mathfrak{A}_2$
will be denoted by $\cL^+(\mathfrak{A}_1,\mathfrak{A}_2)$ ($\cL^+_1(\mathfrak{A}_1,\mathfrak{A}_2)$, respectively).

From now on we make the assumption that $\mathfrak{A}_2$ is equal to $\cB(\cH)$ for some Hilbert space $\cH$.
The important class of positive maps is formed by decomposable maps. They are defined as follows.
Let $T \in \cL^+(\mathfrak{A}, \cB(\cH))$. We say that $T$ is co-positive if $t \circ T$ is completely positive, where $t$ is the transpose map on $\cB(\cH)$. $T$ is decomposable if $T$ is the sum of a completely positive map and a co-positive map. Otherwise, $T$ is called an indecomposable map.

To speak about the structure of  $\cL^+(\qA, \cB(\cH))$ we need some preliminaries.  To this end,
firstly, we note that
in his pioneering work on Banach spaces, Grothendieck \cite{Gro} observed
the links between tensor products and mapping spaces. To describe this, we
will select certain results from the theory of tensor products of Banach
spaces. The point is that the  synthesis  of the linear structure of tensor
products with a topology is not unique - namely, there are many ``good''
cross-norms (cf \cite{Tak}) (the same can be said about the  synthesis  of the linear structure of tensor product and an order, see the next Section). However, among them, there is the projective
norm which gives rise to the projective tensor product and this tensor
product linearizes bounded bilinear mappings just as the algebraic tensor
product linearizes bilinear mappings (see \cite{Ryan}).

Let $X$, $Y$ be Banach algebras.
We denote by $X \odot Y$ the algebraic tensor product of $X$ and $Y$
 (algebraic tensor product of two $^*$-Banach algebras is defined as
tensor product of two vector spaces with $^*$-algebraic structure
determined by the two factors; so the topological questions are not
considered).
We consider the following (projective) norm on $X \odot Y$
\begin{equation}
\pi(u) = \inf \{ \sum_{i=1}^{n}\Vert x_i \Vert \Vert y_i\Vert: \quad u =
\sum_{i=1}^n x_i \otimes y_i \}.
\end{equation}

We denote by $X \otimes_{\pi}Y$ the completion of $X\odot Y$ with respect
to the projective norm $\pi$ and this Banach space will be referred to as
the projective tensor product of the Banach spaces $X$ and $Y$.
Denote by $\mathfrak{B}(X \times Y)$ the Banach space of bounded bilinear
mappings $B$ from $X \times Y$ into the field of scalars with the norm
given by $||B|| = \sup \{ |B(x,y)|; \Vert x \Vert \leq 1, \Vert y \Vert
\leq 1 \}$.
Note (for all details see \cite{Ryan}), that with each bounded bilinear
form $B \in \mathfrak{B}(X \times Y)$ there is an associated operator $L_B
\in \cL(X, Y^*)$ defined by $\langle y, L_B(x)\rangle = B(x,y)$.  The
mapping $B \mapsto L_B$ is an isometric isomorphism between the spaces
$\mathfrak{B}(X \times Y)$ and $\cL(X, Y^*)$. Hence, there is an
identification
\begin{equation}
\label{lala}
(X \otimes_{\pi} Y)^* = \cL(X, Y^*),
\end{equation}
such that the action of an operator $S:X \to Y^*$ as a linear functional
on $X \otimes_{\pi}Y$ is given by
\begin{equation}
\label{haha}
\langle \sum_{i=1}^n x_i \otimes y_i , S \rangle = \sum_{i=1}^n \langle
y_i, Sx_i \rangle.
\end{equation}

 We wish to
complete the presented compilation by
 recalling the another of St{\o}rmer's result (see \cite{St1}) which will be
the crucial in our work. Moreover, it can serve as an illustration on the
given material as well as to indicate that relation (\ref{lala}) is
very relevant to an analysis of positive maps. To present the above
result we need some preparations.

Let $\mathfrak{A}$ be a norm closed self-adjoint subspace of bounded
operators on a Hilbert space $\mathcal K$ containing the identity operator on
$\mathcal K$. $\mathfrak T$ will denote the set of trace class operators
on a Hilbert space $\cH$. Recall that $\cB(\cH) \ni x \mapsto x^t \in \cB(\cH)$ stands
for the transpose map of $\cB(\cH)$ with respect to some orthonormal
basis.
The set of all linear bounded  maps $\phi: \mathfrak{A} \to
\cB(\cH)$ will be denoted by $\cL(\mathfrak{A}, \cB(\cH))$.
 Finally, we denote by ${\mathfrak A} \odot {\mathfrak T}$ the algebraic
tensor product of $\mathfrak A$ and $\mathfrak T$
  where  ${\mathfrak A} {\otimes}_{\pi} \mathfrak T$ means its Banach space
closure under the projective norm defined by

\begin{equation}
\label{4}
\pi(x) = \inf \{ \sum_{i=1}^n \Vert a_i \Vert  \Vert b_i \Vert_1: x =
\sum_{i=1}^n a_i \otimes b_i, \ a_i \in {\mathfrak A}, \ b_i \in
{\mathfrak T} \},
\end{equation}
where $\Vert \cdot \Vert_1$ stands for the trace norm. Now, we can quote
(see \cite{Str})
\begin{lemma}
\label{pierwszy lemat}
There is an isometric isomorphism $\phi \mapsto \tilde{\phi}$ between
$\cL({\mathfrak A}, \cB(\cH))$ and
$({\mathfrak A} {\otimes}_{\pi} {\mathfrak T})^*$ given by
\begin{equation}
\label{5}
(\tilde{\phi})(\sum_{i=1}^n a_i\otimes b_i) = \sum_{i=1}^n
\Tr(\phi(a_i)b^t_i),
\end{equation}
where $\sum_{i=1}^n a_i\otimes b_i \in {\mathfrak A}\odot {\mathfrak T}$.

Furthermore, $ \phi \in \cL^+({\mathfrak A}, \cB(\cH))$ if and only if
$\tilde{\phi}$ is positive on ${\mathfrak A}^+ {\otimes}_{\pi} {\mathfrak
T}^+$.
\end{lemma}

\section{Structure of $\cL^+(\qA, \cB(\cH))$}
We begin with some preliminaries on tensor products of ordered Banach spaces. We already note in the previous section that the synthesis  of the linear structure  of a tensor product with an order is not unique.  Following Wittstock \cite{witt} we start with,
\begin{definition} (\cite{witt})
Let $E$, $F$ be ordered linear spaces with proper cones $E_+$, $F_+$. We call a cone $C_{\alpha} \subset E\otimes F$ a tensor cone and write $E\otimes_{\alpha} F \equiv (E \otimes F, C_{\alpha})$ if the canonical bilinear mappings $\omega : E \times F \to E\otimes_{\alpha}F$ and $\omega^* : E^* \times F^* \to (E\otimes_{\alpha}F)^*$ are positive. Thus, $x\otimes y \in C_{\alpha}$  for all $x\in E_+$, $y \in F_+$; and 
$x^*\otimes y^* \in C_{\alpha}^*$ for all $x^*\in E_+^*$, $y^* \in F_+^*$, where $E_+^*$, $F_+^*$, and $C_{\alpha}^*$ denote the dual cones.
\end{definition}

There are two distinguished cones:
\begin{definition}(\cite{witt})
The projective cone $C_p$:
\begin{equation}
C_p = conv(E_+ \otimes F_+) = \{ \sum_{i=1}^n x_i \otimes y_i; x_i \in E_+, y_i \in F_+, n \in \Nn \}
\end{equation}
\end{definition}
and
\begin{definition}(\cite{witt})
The injective cone $C_i$:
\begin{equation}
C_i = \{ t \in E \otimes F; \left\langle t, E^*_+ \otimes F^*_+\right\rangle  \geq 0 \}.
\end{equation}
\end{definition}
One has, see Proposition 1.14 in \cite{witt}:
\begin{proposition}(\cite{witt})
If $C_{\alpha}$ is a tensor cone, then 
\begin{equation} C_p \subset C_{\alpha} \subset C_i.
\end{equation}
\end{proposition}

Combining  Grothendieck's result,  (\ref{haha}), (\ref{lala}) and Lemma \ref{pierwszy lemat}  with the concept of tensor cones, we get the following selection recipe:
\begin{equation}
\label{3.4}
P_{\alpha} = \{ T \in \cL(\qA, \cB(\cH)); \ \tilde{T}(\sum_{i=1}^n a_i \otimes b_i) = \sum_{i=1}^n \Tr T(a_i)b_i^t \geq 0 \ {\rm for \ any} \  \sum_{i=1}^n a_i \otimes b_i \in C_{\alpha}\},
\end{equation}
 where   $C_{\alpha}$ is a tensor cone in $\qA \otimes_{\pi} \mathfrak{T}$.

\begin{remark}
\label{remark}
In other words, $P_{\alpha}$ is a result of selecting such linear  mappings which have ``nice'' positive functionals, where positivity is defined by  a  selected  tensor cone. It is worth pointing out that there are, in general, many  tensor cones. Consequently, there are many classes of positive maps. Furthermore,  Proposition \ref{3.4} implies:
\begin{equation}
P_i \subset P_{\alpha} \subset P_p.
\end{equation}
\end{remark}
To study the structure of $\cL^+(\qA, \cB(\cH))$ we proceed to a description of some selected cones $C_{\alpha}$ and positive maps in $P_{\alpha}$.
\subsection{{Positive maps determined by $C_p$}}

According to the  recipe given by (\ref{3.4}) we are interested in the following maps:
\begin{equation}
\label{3.5}
P_{p} = \{ T \in \cL(\qA, \cB(\cH)); \ \tilde{T}(\sum_{i=1}^n a_i \otimes b_i) = \Tr T(a_i)b_i^t \geq 0  \ {\rm for \ any} \sum_{i=1}^n a_i \otimes b_i \in C_{p}\}.
\end{equation}
 It is a simple matter to observe that in (\ref{3.5}) a linear bounded map $T$ should satisfy:
\begin{equation}
\Tr  T(a)b \geq 0,
\end{equation}
for any $a \in \qA_+$ and $b \in \mathfrak{T}_+$. But, it is nothing else but the definition of a positive map. Consequently, the smallest tensor cone $C_p$ defines the largest class of positive maps - just the set of all positive maps. We note that this result can be inferred from Stormer's paper, see Lemma \ref{pierwszy lemat} and/or \cite{Str}.

\subsection{Positive maps determined by $C_{cp}$}
We define
\begin{equation} 
C_{cp} = (\qA \otimes_{\pi} \mathfrak{T})_+.
\end{equation}
Firstly, we note that as $\qA$ and $\mathfrak{T}$ are $^*$-Banach algebras, it is easy to check that then $\qA \otimes_{\pi} \mathfrak{T}$ is also $^*$-Banach algebra under the natural involution $(x \otimes y)^* = x^* \otimes y^*$,   cf \cite{Tak},  Section IV.4. Thus $C_{cp}$ is a cone.
Clearly $x \otimes y \in (\qA \otimes \mathfrak{T})_+$ if $x \in \qA_+$ and $y \in \mathfrak{T}_+$.
Furthermore, for $\varphi \in \qA^*_+$ and $\Tr a(\cdot) \in \mathfrak{T}^*_+$, $a \in \cB(\cH)_+$, one has
\begin{equation}
<\varphi \otimes \Tr a(\cdot), \sum_{i,j} x^*_ix_j \otimes \varrho^*_i\varrho_j> =
\sum_{i,j} \varphi(x^*_ix_j) \Tr \varrho^*_i  \varrho_j a \geq 0,
\end{equation}
where $x_i \in \qA$, $\varrho_i \in \mathfrak{T}$, and where we have used the Schur product theorem for the Hadamard product. Thus, $C_{cp}$ is a tensor cone, and we can use the recipe (\ref{3.4}).  So, for $a_i  \in \qA$, $b_i \in \mathfrak{T}$:
\begin{equation}
\label{3.10}
\tilde{T}(\sum_{i,j}^k (a_i \otimes b_i)^*(a_j \otimes b_j))
= \sum_{i,j}^k \Tr T(a^*_ia_j) (b^*_i b_j)^t = \sum_{i,j}^k \sum_n ((b_i)^t e_n, T(a^*_ia_j) (b_j)^te_n)
\end{equation}
Let $(b_j)^t = |f_j><x| \ (\in \mathfrak{T})$, where $f_j \in \cH$, $x \in \cH$ such that $||x|| =1$. Then,
\begin{equation} (\ref{3.10})  = \sum_{i,j}^k \sum_n \overline{(x,e_n)} (f_i, T(a^*_ia_j) f_j) (x,e_n)
= \sum_{i,j}^k (f_i, T(a^*_ia_j) f_j).
\end{equation}
Now, put $f_i = c_i g$, where $c_i \in \cB(\cH)$, $g \in \cH$. Then,
\begin{equation}
\label{3.12}
(\ref{3.10})  = \sum_{i,j}^k (c_i g , T(a^*_ia_j) c_j g) \geq 0
\end{equation}
It follows immediately  from Corollary IV.3.4 in \cite{Tak} that $T$ is a completely positive map. Consequently, the cone $C_{cp}$ is selecting cp maps $P_{cp}$.  Again, we note that a similar  characterization of cp maps could be inferred from St{\o}rmer's paper \cite{Str}.

\subsection{Positive maps determined by $C_i$}
We have seen that the tensor cone $C_i$ is the largest tensor cone. It was already noted by Stinespring \cite{stines} - see the nicely elaborated the Stinespring example in \cite{witt} -  that  in general, $C_i$  is not equal to $C_{cp}$.
Now we wish to describe the corresponding maps which will be denoted by $P_i$. We  note that the condition  defining the cone $C_i$ can be written as

\begin{equation}
\label{3.13}
\sum_{i=1}^n \varphi(a_i) \Tr \varrho_i a \geq 0,
\end{equation}
for any $\varphi \in \qA^*_+$, $a \in \cB(\cH)_+$, where we took $\qA \otimes \mathfrak{T} \ni t = \sum_{i=1}^n a_i \otimes \varrho_i$.

We say that $t = \sum_{i=1}^n a_i \otimes \varrho_i$ satisfying condition (\ref{3.13}) is st-positive (st stands for simple tensor) and denote $t \geq^{st} 0$. We note that the property $t \geq^{st}0$ is, in general, weaker than standard positivity $t \geq 0$, see the example below. Taking into account that  $P_i \subset P_{cp}$ , cf Remark \ref{remark} and the previous subsection,  we infer that the recipe (\ref{3.4}) leads to very regular completely positive maps which will be denoted by 
\begin{equation}
\label{p_i}
P_i  = \{ T:   \sum_i \Tr T(a_i) \varrho_i \geq0\ \ {\rm for} \  \sum_i a_i \otimes \varrho_i \geq^{st} 0 \}. 
\end{equation}

To see  this regularity explicitly as well as the difference between $\geq 0$ and $\geq^{st} 0$, we give an example which can be considered as a continuation of the above mentioned Stinespring's  example.
\begin{example}
Let $\qA$ be equal to $\cB(\cH)$ with $\cH ; \  \dim \cH = n < \infty$.  From the condition (\ref{3.13}) we have
\begin{equation}
\label{3.14}
\sum_i \Tr \varrho_{\varphi} a_i \Tr \varrho_i a \geq 0,
\end{equation}
for any $\varrho_{\varphi} \in \cB(\cH)_+$, $ a \in \cB(\cH)_+$, where $t$  was taken to be of the form $t = \sum_i a_i \otimes \varrho_i$ with $a_i, \varrho_i \in \cB(\cH)$. Put $\varrho_{\varphi} = |f><f|$ and $a = |g><g|$ for $f,g \in \cH$. Then, (\ref{3.14}) leads to
\begin{equation}
\sum_i(f,a_if)(g, \varrho_i g) = \sum_i (f \otimes g, a_i \otimes \varrho_i f \otimes g) \geq 0,
\end{equation}
for any $f,g \in \cH$.  Thus  $\sum_i a_i \otimes \varrho_i \geq^{st} 0$ implies the positivity of $\sum a_i \otimes \varrho_i$ only on simple tensors of $\cH \otimes \cH$. We emphasize that this condition is weaker than the standard positivity $\sum_i a_i \otimes \varrho_i \geq 0$.

To explain the regularity imposed by the cone $C_i$ we recall that  $P_i \subset P_{cp}$, cf Remark \ref{remark}.
In particular, $T \in P_i$ can be written in the form (usually called the Kraus decomposition),
\begin{equation}
T(a) = \sum_k V^*_k a V_k,
\end{equation}
where $V_k \in \cB(\cH)$, see Theorem 4.1.8 in \cite{St}. 
As a first observation, we consider the case: $V_k = |f_k><g_k|$ with $f_k, g_k \in \cH$.  Then,
\begin{equation}
\sum_{i,k} \Tr V^*_k a_i V_k \varrho_i = \sum_{i,k,l} (e_l, V^*_k a_i V_k\varrho_i e_l) = \sum_{i,k,l} (|f_k><g_k|e_l, a_i 
|f_k><g_k| \varrho_i e_l) 
\end{equation}
$$ = \sum_{i,k,l} (f_k,  a_if_k)(\varrho_i^*g_k,e_l)(e_l,g_k) = \sum_{i,k} (f_k, a_i f_k)(g_k, \varrho_i g_k)$$
$$ = \sum_k (f_k\otimes g_k, \left( \sum_i a_i \otimes \varrho_i \right) f_k \otimes g_k).$$
Consequently, $T(\cdot) = \sum_i V_k^* (\cdot )V_k$ with $V_k = |f_k><g_k|$ is a regular cp map which is in $P_i$.
It is worth pointing out that one-rank operators $V_k$ ensures applicability of st positivity and that such maps $T(\cdot) = \sum_i V_k^* (\cdot )V_k$ are sometimes called super positive maps, cf  Definition 5.1.2 in \cite{St}.  Now, we turn to the case $V_k = \sum_m^M|f^k_m><g^k_m|$ with $M \geq2$.
\begin{equation}
\label{rownanie}
\sum_{k,l,i}(e_l, V^*_k a_iV_k\varrho_i e_l) = \sum_{i,k,l} (\sum_m|f^k_m><g^k_m| e_l, a_i \sum_n |f^k_n><g^k_n| \varrho_i e_l)
\end{equation}
$$= \sum_{i,k,l,m,n} (e_l,g^k_m)(f^k_m, a_i f^k_n)(g_n^k,\varrho_i e_l)=
\sum_{k,m,n} (f^k_m\otimes g^k_n, \left( \sum_i a_i \otimes\varrho_i \right)  f^k_n\otimes g^k_m).$$
We see at once that st. positivity of $\sum_i a_i \otimes \varrho_i$ does not ensure the positivity of (\ref{rownanie}). Consequently, such cp maps are not, in general, in $P_i$.
\end{example}
\subsection{Positive maps determined by $C_{cp} \cap id\otimes t(C_{cp})$}
We begin by an examination of the set $C_d \equiv C_{cp} \cap id\otimes t(C_{cp})$, where $t$ as before, stands for the transposition. We first note that $C_d$ is a cone. To deduce that $C_d$ is a tensor cone , we note that $x\otimes y \in C_d$ for all $x \in \qA_+$ and $y \in \mathfrak{T}_+$. Subsequently, we observe that  one has $x^* \otimes y^* \in C^*_d$ for all  $x^* \in \qA^*_+$ and $y^*  \in \mathfrak{T}^*_+$. Thus, $C_d$ is a tensor cone. Clearly, 
\begin{equation}
\label{jeden}
C_p \subset C_d \subset C_{cp} \quad \rm{and} \quad P_{cp} \subset P_d \subset  P_p,
\end{equation}
where $P_d$ stands for those positive maps which are determined by the cone $C_d$, i.e.
\begin{equation}
P_d = \{ T: \sum_{i,j =1}^n  \Tr T(a^*_ia_j )(b^*_i b_j)^t \geq 0\ {\rm{and}} \ \sum_{i,j =1}^n  \Tr T(a^*_ia_j )(b^*_i b_j)\geq 0,
\end{equation}
where $a_i \in \qA$, and $b_i \in \mathfrak{T}$.

Further, let us consider the set $C_{ccp} \equiv id \otimes t(C_{cp})$. It follows by similar arguments as those  employed in Subsection  3.2
that $C_{ccp}$ is a tensor come. Hence, it is easy to check that
\begin{equation}
\label{dwa}
C_p \subset C_d \subset C_{ccp} \quad \rm{and} \quad P_{ccp} \subset P_d \subset  P_p,
\end{equation}
where $P_{ccp}$ stands for the set of all co-positive maps.
Consequently, (\ref{dwa}) and (\ref{jeden})  lead to 
\begin{equation}
P_d \supseteq P_{cp} \cup P_{ccp}.
\end{equation}
In other words, the cone $C_d$ determines decomposable maps $P_d$.

To finish this subsection we have to examine the question whether $C_d$ is always non-trivial, i.e. whether the inclusion $C_p \subset C_d$ is the proper one. To answer this question we give:
\begin{example}
\label{przyklad1}
We assume that $\qA = \cB(\cH)$ and that $\dim \cH \leq 3$. To study the non-triviality of $(\cB(\cH) \otimes \cB(H))_+ \cap
id \otimes t((\cB(\cH) \otimes \cB(H))_+)$ it is enough to examine one dimensional (orthogonal) projectors of the form $|f><f|$ with $f \in \cH \otimes \cH$. To see this we begin with two observations: 
\begin{observation}
For $\cH \otimes \cH \ni h = \sum_i f_i \otimes g_i$ one has
\begin{equation}
|h><h| = | \sum_i f_i \otimes g_i><\sum_j f_j\otimes g_j| = \sum_{i,j}|f_i><f_j| \otimes |g_i><g_j|.
\end{equation}
\end{observation}
\begin{observation}
Any $h \in \cH \otimes \cH$ can be written in the form $h = \sum_i v_i \otimes e_i  \ (= \sum e_i \otimes w_i )$ where $v_i, w_i \in \cH$ and $\{e_i \}$ is a basis in $\cH$.
\end{observation}
Proofs of the above observations are easy  and they are left to the reader.

Let $|f><f|$ be given.  Then, the above  Observations and (\ref{3.4}) imply
\begin{equation}
\tilde{T}(|f><f|) = \sum_{k,l} \Tr T(|e_k><e_l|) |f_l><f_k|,
\end{equation}
where $f= e_k \otimes f_k$, $\{ e_i \}$ a basis in $\cH$. 
To see that that the Woronowicz scheme for description of positive maps of low dimensional matrix algebras  is reproduced, \cite{Wor},  we note,  in the ``matrix terms'',  that
\begin{equation}
\left\lbrace \sum_l T(|e_k><e_l|) |f_l><f_k| \right\rbrace _{kk} = \sum_l \left\lbrace T(|e_k><e_l|) \right\rbrace _{kl}
\left\lbrace |f_l><f_k|  \right\rbrace_{lk}, 
\end{equation}
and
$$ \Tr_{\cH\oplus \cH} \sum_l [T(|e_k><e_l|)] Q_{lm} = \sum_{k,l} \Tr_{\cH} T(|e_k><e_l|) \cdot |f_l><f_k|,$$
where, using Woronowicz's notation, $Q$ stands for the operator $\{ |f_k><f_l|\}_{kl}$.
Consequently, to infer the triviality of the cone $C_d$ for dimension $2$ from the Woronowicz result, Theorem 1.1 and Section 2 in \cite{Wor},
it is enough to reproduce, in the considered context, Woronowicz's argument leading to ``simple vectors''. To this end, we proceed to show that,  for $\cH; \dim \cH =2$,  
 $id \otimes t(h><h|)$ is positive only if $h=f\otimes g$, $f,g \in \cH$, so when $|h><h| \in \cB(\cH)_+ \otimes \cB(\cH)_+$.  Let $x = \sum_m w_m \otimes e_m  \in \cH \otimes \cH$, and $h = \sum_i v_i \otimes e_i \in \cH \otimes \cH$. Then, 
 \begin{equation}
 (x, id \otimes t (h><h|) x ) = \sum_{i,j,m,n} (w_n \otimes e_n, |v_i><v_j| \otimes |e_j><e_i| w_m \otimes e_m)
 \end{equation}
 $$=\sum_{i,j,m,n} (w_n,v_i)(e_n,e_j)(v_j,w_m)(e_i,e_m) =\sum_{i,j} (w_j,v_i)(v_j,w_i),$$
 where both observations were used.
 
 For, $\cH; \dim \cH = 2$ one gets
\begin{equation}
\label{Q1}
(w_1,v_1)(v_1,w_1) + ( w_2,v_1)(v_2,w_1) + (w_1,v_2)(v_1,w_2) + (w_2,v_2)(v_2,w_2) \geq 0,
\end{equation}
while for $\cH; \dim \cH = 3$
\begin{equation}
\label{Q2}
(w_1,v_1)(v_1,w_1) + ( w_2,v_1)(v_2,w_1) + (w_1,v_2)(v_1,w_2) + (w_2,v_2)(v_2,w_2) + (w_1,v_3)(v_1,w_3)
\end{equation}
$$ + (w_2,v_3)(v_2,w_3) + (w_3,f_3)(v_3,w_3) + (w_3,v_1)(v_3,w_1) + (w_3,v_2)(v_3,w_2)\geq0.
$$
 We are looking for $\{v_i \}$ such that (\ref{Q1}) ( respectively (\ref{Q2}) ) are satisfied for arbitrary $\{w_i\}$.
 
 Let's consider (\ref{Q1}) in detail. As $\{w_i\}$ was arbitrary we may change $w_1$ for $\lambda w_1$, $\lambda \in \Rn$
 to arrive at
 \begin{equation}
 \label{Q3}
 \forall_{\lambda}  \quad F(\lambda) = \lambda^2 a + \lambda b +c \geq 0,
 \end{equation}
 where $a = (w_1,v_1)(v_1, w_1)$, $b= (w_2,v_1)(v_2,w_1) + \overline{(w_2,v_1)(v_2,w_1)}$, $c = (w_2,v_2)(v_2,w_2)$.
 (\ref{Q3}) holds if $c \leq \frac{b^2}{4a}$. As $\{ w_i \}$ are arbitrary, we may take $w_2 \bot v_1$.  Then, $b=0$. So $c=0$. As $w_1$ is arbitrary, then $v_1 = \lambda v_2$, where $\lambda \in \Cn$.
 However, this implies that $h$ is a simple tensor.

We now turn to the case $\dim \cH = 3$. It is easily seen that the above arguments fail when we have the 3 dimensional case. 
Fortunately, Choi gave an explicit construction of a matrix $U \in (M_3(\Cn) \otimes M_3(\Cn))^+$  such that $U\neq M_3(\Cn)^+ \otimes M_3(\Cn)^+,$ see \cite{Choi2}; another example was given in \cite{Ster}.
 So, it is enough to see that the cone $C_d$ is not trivial for $\cH; \dim \cH = 3$.
 \end{example}
 \subsection{k-positive maps}
 To complete this section we wish to examine $k$-positive maps.  Our first observation is that we can apply the strategy described in subsection $3.1$ with the following modifications:
\begin{equation} 
 \qA \rightarrow M_n(\qA)\cong M_n(\Cn)\otimes \qA,
 \end{equation}
 and
 \begin{equation}
 B(\cH)  \rightarrow M_n(B(\cH))\cong M_n(\Cn) \otimes B(\cH).
 \end{equation}
 Consequently, we will examine
 \begin{equation}
 \cL\left( M_n(\qA), M_n(B(\cH)) \right) \cong\left( M_n(\qA) \otimes_{\pi} M_n(B(\cH))_*\right) ^*.
 \end{equation}
 To proceed with the analysis of $k$-positive maps we remind the reader that any $A \in M_n(\qA)$ has the unique representation:
 \begin{equation}
 \label{3.34}
 A = \sum_{i,j=1}^n  \epsilon_{ij} \otimes a_{ij},
 \end{equation}
 where $a_{ij} \in \qA$, and $\{ \epsilon_{ij} \}_{i,j =1}^n$ is a matrix unit, ie $(\epsilon_{ij})^* = \epsilon_{ji}$, $\epsilon_{ij}\epsilon_{kl} = \delta_{jk}\epsilon_{il}$, $\sum_{i=1}^n \epsilon_{ii} = 1$ cf \cite{Wegge}, Appendix $T$.
 The above representation of $A \in M_n(\qA)$ combined with Kaplan's arguments, see Proposition 1.1 in \cite{kaplan}, plus obvious modifications,  give an order isomorphism $\mathfrak{I}: M_n(B(\cH)_*) \to M_n(B(\cH))_*$ given by 
 \begin{equation}
 \left( \mathfrak{I}([\varrho_{ij}])\right) = \sum_{ij=1}^n \varrho_{ij}(a_{ij}),
 \end{equation}
 where $[\varrho_{ij}] \in M_n(B(\cH)_*)$, $[a_{ij}] \in M_n(B(\cH))$. Then, using the recipe analogous to that given by \ref{3.4} one has the following prescription for selection of $k$-positive maps:
 \begin{eqnarray*}
 (*) \quad \{   T_k \equiv id_{M_k(\Cn)} \otimes T \in \cL(M_k(\qA), M_k(B(\cH))); \tilde{T_k}(\sum_{p=1}^{p_0} [a^p_{ij}] \otimes [\varrho^p_{ij}]) \\
 = \sum_{p=1}^{p_0} \sum_{ij}^k \Tr T(a^p_{ij}) \varrho^{p,t}_{ij} \geq 0,
 \rm{for \ any } \  \sum_{p=1}^{p_0} [a^p_{ij}] \otimes [\varrho^p_{ij}] \ {\rm{in}} \  M_k(\qA)^+ \otimes M_k(B(\cH)_*)^+ \} .
 \end{eqnarray*}
 
 To verify the above procedure, $(*)$,  it is enough to note that $[\varrho_{ij}] \in M_k(B(\cH)_*)^+= B(\Cn^n \otimes \cH)_*^+$ iff $\varrho_{ij} = \sum_{l=1}^k \sigma_{li}^*\sigma_{lj}$, where the unique representation   (\ref{3.34}) together with properties of trace class and Hilbert-Schmidt operators  were  used.  Hence
 \begin{equation}
 \label{36}
 0\leq \sum_{i,j =1}^k \Tr \sum_{m=1}^kT(a^p_{ij}) \sigma^t_{mj} (\sigma^t_{mi})^* = \sum _{ij = 1^k} \Tr \sum_m (\sigma^t_{mi})^* T(a^p_{ij}) \sigma_{mj}^t.
 \end{equation}
 So $[T(a^p_{ij})]_{ij=1}^k \geq 0$ and $T$ is $k$-positive.
 
 Finally we note that  dropping the restriction for the upper bound of summation in (\ref{36}) one gets
 \begin{equation}
 \label{37}
 0\leq \sum_{i,j,m =1}^l \Tr T(a_{ij}) \sigma^t_{mj} (\sigma^t_{mi})^* = \sum _{ijm} \Tr (\sigma^t_{mi})^* T(a_{ij}) \sigma_{mj},
 \end{equation}
 where $l$ is an arbitrary natural number. It is easily seen that this is exactly the condition imposed by the cone $C_{cp}$. In other words, this observation  sheds  some new light on the origin of the definition of CP-maps.

\section{Conclusions and remarks}
The principal significance of the Example \ref{przyklad1} is that it clarifies the appearance of non-decomposable maps, i.e. 
if $C_d \neq C_p$ then there is a room for non-decomposable maps. We have seen that both cones $C_d$ and $C_p$ are equal to each other  for $\cH; \dim \cH =2$.

The next important point to note here is that when studying the structure of the set of positive maps we restrict ourselves to a 
few special  cones. In general, there could be other tensor cones. Thus, in principle, there could be other classes of interesting positive maps.

In \cite{WAM_o} the structure of the set of positive maps was examined using the concept of elementary maps. Although this concept was to some extent vague,
it was indicated that  general structural properties should play  a significant role. Here, we have seen that the synthesis  of Grothendieck's idea with 
the order on tensor product is playing a crucial role. 

It is worth pointing out that the presented scheme is  offering another approach to the concept of elementary maps.
To see this, from now on we make the assumption: $\qA = \cB(\cH)$ with $\cH; \dim \cH =n <\infty$, cf \cite{WAM_o}.  We emphasize  that this case is essential for Quantum Information Theory.

We first  note that in the recipe (\ref{3.4}) one can restrict oneself to extremal functionals. As an extreme point in a convex set can be treated as an elementary constituent, we are getting another concept of an elementary map. Define
\begin{equation}
P_{\alpha}^e = \{ T \in \cL(\cB(\cH), \cB(\cH)); \ \tilde{T}(\sum_{i=1}^n a_i \otimes b_i) = \sum_{i=1}^n \Tr T(a_i)b_i^t \geq 0 \}
\end{equation}
for any $\sum_{i=1}^n a_i \otimes b_i \in C_{\alpha}$, where  $\sum_{i=1}^n a_i \otimes b_i  \in \cB(\cH) \odot \cB(\cH)$, $C_{\alpha}$ is a tensor cone in $\cB(\cH) \otimes_{\pi} \cB(\cH)$, and $\tilde{T}$ is assumed to be an extremal functional.
Consequently, elementary maps can be written  as: 
\begin{equation}
\label{4.2}
P_{\alpha}^e = \{ T \in \cL(\cB(\cH), \cB(\cH)); \ (h, (\sum_{i=1}^n a_i \otimes b_i) h) = \sum_{i=1}^n \Tr T(a_i)b_i^t \geq 0 \}
\end{equation}
for any $\sum_{i=1}^n a_i \otimes b_i \in C_{\alpha}$, where  $\sum_{i=1}^n a_i \otimes b_i  \in \cB(\cH) \odot \cB(\cH)$, $C_{\alpha}$ is a tensor cone in $\cB(\cH) \otimes_{\pi} \cB(\cH)$ and $h \in \cH \otimes \cH$.

The interest of this remark is that it provides a recipe for constructing ``elementary'' maps with specified positivity with respect to the selected cone $C_{\alpha}$. Furthermore, in an application  of the above scheme to $2$-positive maps one recoves the concept of \textit{atomic maps}.

The next important point to note here is that the relation between  a cp-map $T$ and the asociated  $C_{cp}$-positive functional $\tilde{T}$ can  be considered as the starting point for a generalization of  Krauss (Stinespring) decomposition, for the corresponding order structures see Section 1.4 in \cite{Arv}.

 To see this we recall  that $\left( M_n(\Cn)\otimes_{\pi} M_n(\Cn)_*) , C_{cp} \right) $ is an involutive Banach algebra and the cone $C_{cp}$  can be used to define the concept of states, cf Chapter I in \cite{Tak}. Therefore, there is a posibility to speak about GNS construction.
Then employing the decomposition theory,  as it was given in Chapter 4, in \cite{BR}, see Section 4.1.1,  one can decompose the $C_{cp}$-positive functional $\tilde{T}$, so also the corresponding map $T$.  In other words, we are getting a generalized decomposition of a cp-map.

Our next remark is clarifying a general form of a positive map (for finite dimensional case). We note that for  any self-adjoint (hermitan)  functional $\tilde{T}$ on a $C^*$-algebra $\qA$ there is a \textit{unique} pair  $\tilde{T}_+$  and $\tilde{T}_-$  of positive functionals such that $\tilde{T}  = \tilde{T}_+  - \tilde{T}_-$ , cf Theorem 3.2.5 in \cite{Ped}. We can use this for  unique decomposition of  a self-adjoint functional  on  $M_n(\Cn)\otimes M_n(\Cn)$. Subsequently, one can  transfer this  decomposition  to $M_n(\Cn)\otimes_{\pi} M_n(\Cn)$.  Thus, we arrive at: a $C_p$-positive functional $\tilde{T}$ has \textit{unique} decomposition $\tilde{T}  = \tilde{T}_+  - \tilde{T}_-$ , where $\tilde{T}_+$, $\tilde{T}_-$ are $C_{cp}$-positive functionals.  

Consequently we got, using our scheme, the following decomposition of a positive map $T$; $T = T_+ - T_-$, where $T_+$, $T_-$ are cp-maps.  To elaborate and clarify the decomposition formula  of $T$  we need some preliminaries, cf. \cite{hou}.

Let $a_1,...,a_k$ ,and $c_1,...,c_i$ be in $\cB(\cH,\cK)$ ($\cH, \cK$ stand for Hilbert spaces). If for each $x \in \cH$, there exists an $l\times k$ complex matrix $\{(\alpha_{i,j}(x) )\}$ such that
\begin{equation}
c_ix = \sum_{j=1}^k \alpha_{i,j}(x) a_j x, \quad i= 1,...,l,
\end{equation}
$(c_1,...,c_l)$ is said to be a locally linear combination of $(a_1,...,a_k)$. If coefficients $\{\alpha_{i,j}(x)\}$ can be taken in such way that the norm $\| \alpha_{i,j}(x)\| \leq 1$, for every $x$, then $(c_1,...,c_l)$ is said to be a contractive linear combination of $(a_1,...,a_k)$.   Employing Hou's result  (see Corollary 2.6 in \cite{hou}) we arrive at
\begin{proposition}
\label{4.1}
A positive map  $T: M_n(\Cn)\otimes M_n(\Cn)$  can be uniquely written in the form $T = T_+ - T_-$ where $T_+$, $T_-$ are cp-maps. In particular $T_+(a) = \sum_r c_rac_r^*$ where $c_r \in M_n(\Cn)$.  Furthermore, there exists $\{d_p\} \subset M_n(\Cn)$ such that
 $(d_1,...,d_l)$ is a contractive locally linear combination of $(c_1,...,c_k)$ and
\begin{equation}
T(a) = T_+(a) - T_-(a)  = \sum_r  c_i v c_i^* - \sum_p d_j v d_j^*,
\end{equation}
for all $a \in M_n(\Cn)$.
\end{proposition}

Finally we note that in \cite{MMO}  the  identification of the injective cone $C_i$ with the cone $C_{cp}$ was incorrect.

\section{Acknowledgments}
The author would like to express his thanks to Louis E. Labuschagne and Marcin Marciniak for several helpful comments. Furthermore, he
 wishes to thank the University of Gdansk where a part of the research was done.


\begin{thebibliography}{99}
\bibitem{paulsen} V. Paulsen, \textit{ Completely bounded maps and operator algebras}, Cambridge University Press, 2002
\bibitem{St} E. St{\o}rmer,\textit{Positive linear maps on operator algebras}, Springer-Verlag, 2013
\bibitem{St1} E. St{\o}rmer, ``Positive linear maps of operator algebras'', \textit{Acta Mathematica} \textbf{110}, 233 (1963)
\bibitem{Choi} M-D. Choi, ``Positive semidefinite biquadratic forms'', \textit{Linear Algebra and Appl.} \textbf{12}, 95 (1975)
\bibitem{Choi1} M-D. Choi, ``Some asorted inequalities for positive linear maps on $C^*$-algebras'',\textit{J. Operator Theory} \textbf{4}, 271 (1980)
\bibitem{Wor} S. L. Woronowicz, ``Positive  maps of low dimensional matrix algebras'', \textit{Rep. Math. Phys.} \textbf{10}, 165 (1976)
\bibitem{WAM}  W. A. Majewski, Transformations between quantum states, \textit{Rep. Math. Phys.}, \textbf{8}, 295 (1975)
\bibitem{WAM2} W. A. Majewski, M. Marciniak, ``On a characterization of a positive maps'', \textit{J.Phys. A, Math. Gen.} \textbf{34} 5863 (2001)
\bibitem{LMM} L. E. Labuschagne, W. A. Majewski, M. Marciniak, ``On k-decoposability of positive maps'', \textit{Expo.Math.} \textbf{24} 103 (2006)
\bibitem{WAM3} W. A. Majewski, ``On the structure of positive maps: Finite-dimensional case'', \textit{J. Math. Phys.} \textbf{53}, 023515 (2012)
\bibitem{MMO} W. A. Majewski, T. Matsuoka, M. Ohya, Characterization of partial positive states and measures of entaglement, \textit{J. Math. Phys.} \textbf{50}, 113509 (2009)
\bibitem{WAM4}  W. A. Majewski, ``On positive maps in quantum information'', \textit{Ros. J. Math. Phys.} \textbf{21}, 362 (2014)
\bibitem{WAM_o} W. A. Majewski, On the origin of non-decomposable maps, arXiv 1706.07945v2 [mathOA]
\bibitem{Gro}A. Grothendieck, ``Products tensoriels topologiques et espaces nuclearies'', \textit{Memoirs of the American Mathematical Society}, vol. \textbf{16} Providence, RI, 1955
\bibitem{Tak}  M. Takesaki, \textit{Theory of Operators Algebras I}, Springer Verlag, 1979
\bibitem{Ryan} R. A. Ryan, \textit{Introduction to Tensor Products of Banach Spaces}, Springer-Verlag, 2002
\bibitem{Str} E. St{\o}rmer, Extension of positive maps into $\cB(\cH)$, \textit{J. Funct. Anal.} \textbf{66}, 235-254 (1986)
\bibitem{witt} G. Wittstock, \textit{Ordered Normed Tensor products}, in \textit{Foundations of Quantum Mechanics and Ordered Linear Spaces}, Springer-Verlag, Lecture Notes in Physics, vol. 29, 1974
\bibitem{stines} W. F. Stinespring, ``Positive functions on $C^*$-algebras'',  \textit{Proc. Am. Math. Soc.} \textbf{6} 211-216 (1955)
\bibitem{Choi2}  M-D. Choi, ``Positive linear maps'', in \textit{Proceedings of Symposia in Pure Mathematics}, vol. \textbf{38}, Part 2, 583-590 (1982)
\bibitem{Ster} E. St{\o}rmer, ``Decomposable positive maps on $C^*$-algebras'', \textit{ Proc. Amer. Math. Soc.} vol. \textbf{86}, 402-404 (1982)
\bibitem{Wegge} N. E. Wegge-Olsen, \textit{$K$-theory and $C^*$-algebras. A friedly approch.} Oxford University Press, 1993
\bibitem{kaplan} A. Kaplan, ``Multi-states on $C^*$-algebra'', \textit{Proc. Amer. Math. Soc.} vol. \textbf{106}, 437-446 (1989)
\bibitem{Arv} W. B. Arveson, Subalgebras of $C^*$-algebras, \textit{Acta Mathematica} \textbf{123}, 142-224 (1970) 
\bibitem{BR} O. Bratteli, D. W. Robinson, \textit{Operator algebras and Quantum Statistical Mechanics I}, Springer-Verlag, 1979.
\bibitem{Ped}  G.K. Pedersen, \textit{$C^*$-algebras and their automorphism groups}, Academic Press, 1979
\bibitem{hou} Jin-Chuan Hou, ``A characterization of positive elementary operators'', \textit{J.Operator Theory} \textbf{39}, 43-58 (1998)
\end{thebibliography}
\end{document}